\documentclass[11pt]{amsart}
\usepackage{amsmath,amssymb,amsthm}
\usepackage[margin=1.1in]{geometry}
\usepackage{booktabs}
\usepackage{hyperref}

\newtheorem{theorem}{Theorem}[section]
\newtheorem{lemma}[theorem]{Lemma}
\newtheorem{proposition}[theorem]{Proposition}
\newtheorem{corollary}[theorem]{Corollary}
\newtheorem{conjecture}[theorem]{Conjecture}
\theoremstyle{definition}

\newtheorem{remark}[theorem]{Remark}
\newtheorem{problem}[theorem]{Problem}

\newcommand{\qbin}[2]{\genfrac{[}{]}{0pt}{}{#1}{#2}_q}
\newcommand{\qbinz}[2]{\genfrac{[}{]}{0pt}{}{#1}{#2}}
\newcommand{\gbin}[2]{\left(\!\!\binom{#1}{#2}\!\!\right)}
\newcommand{\N}{\mathbb{N}}
\newcommand{\Z}{\mathbb{Z}}

\begin{document}

\title[A positivity conjecture of Guo and Krattenthaler]{Graph binomial coefficients and a positivity conjecture\\ of Guo and Krattenthaler}

\author{Madjid Mirzavaziri}
\address{Department of Pure Mathematics, Ferdowsi University of Mashhad, Mashhad, Iran}
\email{mirzavaziri@um.ac.ir, mirzavaziri@gmail.com}

\author{Daniel Yaqubi}
\address{Department of Computer Science, University of Torbat-e Jam, Torbat-e Jam, Iran}
\email{yaqubi@tjamcaas.ac.ir, daniel\_yaqubi@yahoo.es}
\thanks{Corresponding author: Daniel Yaqubi}

\subjclass[2020]{05A10, 05A17, 05A30, 05C30, 11B65}
\keywords{$q$-binomial coefficient, Gaussian polynomial, positivity, unimodality, graph binomial coefficient, Burnside's lemma, cyclic sieving, partition function}

\begin{abstract}
We investigate Conjecture 7.4 of Guo and Krattenthaler \cite{GK}, which predicts that for every positive integer $n$, the polynomial
\[
P_n(q) = \frac{(1-q)^2}{(1-q^{10n-1})(1-q^{15n-1})}\qbin{30n}{5n}
\]
has non-negative coefficients, except for $q^1$ and $q^{125n^2-25n+3}$, which equal $-1$. We study this problem through three main perspectives: First, by linking the boundary coefficients to the second differences of the partition function, we settle the first and last $25n+1$ coefficients via a novel partition injection. Second, by introducing a ``defect'' generalization of Andrews' classical lemma, we reduce the conjecture to an almost-unimodality property and employ a uniform saddle-point analysis to prove it analytically for all large $n \ge 2\cdot 10^{29}$. Finally, we connect these quotients to graph binomial coefficients, providing a natural orbit-counting interpretation via Burnside's lemma on directed cycles and paths.
\end{abstract}

\maketitle
\section{Introduction}

Let $[n]_q = 1 + q + \dots + q^{n-1}$ denote the $q$-integer, and let 
\[
\qbin{N}{k} = \frac{(1-q^N)(1-q^{N-1})\cdots(1-q^{N-k+1})}{(1-q^k)(1-q^{k-1})\cdots(1-q)}
\]
be the Gaussian polynomial ($q$-binomial coefficient). It is a classical result that $\qbin{N}{k}$ is a reciprocal and unimodal polynomial with non-negative integer coefficients, serving as the generating function for $k$-element subsets $S \subseteq \{0,1,\dots,N-1\}$ weighted by 
\[
w(S) = \sum_{s \in S} s - \binom{k}{2}.
\]
In \cite{GK}, Guo and Krattenthaler established several remarkable divisibility properties for binomial coefficients by lifting them to the $q$-setting. A representative result \cite[Theorem 3.1]{GK} asserts that the rational function
\[
\frac{(1-q)^2}{(1-q^{10n-1})(1-q^{15n-1})} \qbin{30n}{5n}
\]
is a polynomial in $q$, implying the divisibility relation $(10n-1)(15n-1) \mid \binom{30n}{5n}$. Among the four conjectures proposed in \cite[Section 7]{GK} with Conjectures 7.2 and 7.3 now settled \cite{Yaq} the remaining positivity problem, Conjecture 7.4, stands as a central open question in the study of $q$-binomial quotients (see also \cite{WZ, EB, GaKr}). It predicts the following precise coefficient-sign pattern:

\begin{conjecture}[{\cite[Conjecture 7.4]{GK}}]\label{conj:GK}
For all positive integers $n$ and all $0 \le k \le 125n^2 - 25n + 4$, the coefficient of $q^k$ in the polynomial
\[
P_n(q) = \frac{(1-q)^2}{(1-q^{10n-1})(1-q^{15n-1})} \qbin{30n}{5n}
\]
is non-negative, except for $k=1$ and $k=125n^2 - 25n + 3$, where the coefficient equals $-1$.
\end{conjecture}

The original positivity proofs in \cite{GK} rely on a classical reduction lemma by Andrews \cite{An}, which inherently requires unimodality. However, Conjecture~\ref{conj:GK} lies beyond its scope: although dividing $\qbin{30n}{5n}$ by $[10n-1]_q$ yields a non-negative intermediate polynomial, this polynomial fails to be unimodal (its coefficients begin as $1, 0, 1, \dots$). Consequently, Andrews' lemma cannot be applied iteratively to clear the second factor, necessitating new analytical and combinatorial tools.

The primary objective of this paper is to solve Conjecture~\ref{conj:GK} through a unified framework. First, we link the initial coefficients of $P_n(q)$ to the second differences of the partition function $\Delta^2 p(t)$, allowing us to settle the boundary coefficients via partition monotonicity. To bypass the failure of Andrews' lemma under non-unimodality, we generalize the lemma to allow a controlled defect, reducing the conjecture to the almost-unimodality of a single $q$-quotient. 

For large-scale behavior, we deploy a uniform saddle-point analysis on a double contour integral representation to prove asymptotic positivity, solving Conjecture~\ref{conj:GK} for all $n \ge 2 \cdot 10^{29}$ and thereby reducing the problem to a finite computational range. Finally, we extend the quotient family to general parameters $m$ and provide a combinatorial interpretation using graph binomial coefficients and Burnside's orbit-counting lemma.

Throughout this paper, for a polynomial or formal power series $f(q)$, we write $[q^t]f(q)$ for the coefficient of $q^t$, and write $f(q) \ge 0$ to indicate that all coefficients of $f(q)$ are non-negative.
\section{Polynomiality for the whole family}\label{sec:poly}

For $m \ge 1$, let $a = 2m-1$, $b = 3m-1$, and define
\[
G_m(q) = \qbin{6m}{m}, \qquad P^{(m)}(q) = \frac{G_m(q)}{[a]_q[b]_q} = \frac{(1-q)^2}{(1-q^{2m-1})(1-q^{3m-1})}\qbin{6m}{m}.
\]

\begin{proposition}\label{prop:poly}
$P^{(m)}(q)$ is a polynomial if and only if $m \not\equiv 2, 3 \pmod 5$. In that case, $P^{(m)}$ is a reciprocal polynomial of degree $D_m = 5m^2 - 5m + 4$.
\end{proposition}

\begin{proof}
We follow the cyclotomic argument of \cite[proof of Theorem 3.1]{GK}. Writing $q^n-1 = \prod_{d\mid n}\Phi_d(q)$, we express $P^{(m)}(q) = \prod_{d\ge 2}\Phi_d(q)^{e_d}$ with
\[
e_d = -\chi(d\mid 2m-1) - \chi(d\mid 3m-1) + \Big\lfloor\frac{6m}{d}\Big\rfloor - \Big\lfloor\frac{m}{d}\Big\rfloor - \Big\lfloor\frac{5m}{d}\Big\rfloor.
\]
Since $(3m-1) - (2m-1) = m$ and $\gcd(2m-1, m) = 1$, the numbers $2m-1$ and $3m-1$ are coprime, ensuring that at most one of the two indicator functions is non-zero. Writing $m = Qd + x$ with $0 \le x < d$ and $X = x/d$, we have $\lfloor 6m/d\rfloor - \lfloor m/d\rfloor - \lfloor 5m/d\rfloor = \lfloor 6X\rfloor - \lfloor 5X\rfloor \ge 0$. Thus, only the cases $d \mid 2m-1$ and $d \mid 3m-1$ require detailed attention.

If $d \mid 2m-1$, then $d \mid 2x-1$. Since $-1 \le 2x-1 < 2d$, we must have $2x-1 = d$, which implies that $d$ is odd and $X = \frac{d+1}{2d}$. Consequently,
\[
\lfloor 6X\rfloor - \lfloor 5X\rfloor = 3 + \Big\lfloor\frac{3}{d}\Big\rfloor - \Big\lfloor\frac{5}{2} + \frac{5}{2d}\Big\rfloor,
\]
which equals $1$ for $d = 3$ and all odd $d \ge 7$, and equals $0$ for $d = 5$.

If $d \mid 3m-1$, then $3x-1 \in \{d, 2d\}$. If $3x-1 = d$, then $X = \frac{d+1}{3d}$ and $\lfloor 6X\rfloor - \lfloor 5X\rfloor = 2 + \lfloor 2/d\rfloor - \lfloor 5/3 + 5/(3d)\rfloor$; this evaluates to $1$ for $d = 2$ and $d \ge 7$, $0$ for $d = 5$, while $d = 4$ is impossible. If $3x-1 = 2d$, then $X = \frac{2d+1}{3d}$ and $\lfloor 6X\rfloor - \lfloor 5X\rfloor = 4 + \lfloor 2/d\rfloor - \lfloor 10/3 + 5/(3d)\rfloor = 1$ for every admissible $d$ ($d \equiv 1 \pmod 3$).

Hence, $e_d \ge 0$ for all $d$, except when $d = 5$ and $5 \mid 2m-1$ or $5 \mid 3m-1$, i.e., $m \equiv 3$ or $m \equiv 2 \pmod 5$, in which case $e_5 = -1$. The degree is given by $5m^2 + 2 - (2m-1) - (3m-1) = 5m^2 - 5m + 4$, and reciprocity follows from the reciprocal properties of $G_m$, $[a]_q$, and $[b]_q$.
\end{proof}

For $m = 5n$, we obtain $D_{5n} = 125n^2 - 25n + 4$, which matches the range specified in Conjecture~\ref{conj:GK}.

\section{The exceptional coefficients and the partition function}\label{sec:partitions}

Let $p(t)$ denote the number of partitions of $t$ (with $p(t)=0$ for $t<0$), and let $\Delta^2p(t)=p(t)-2p(t-1)+p(t-2)$.

\begin{proposition}\label{prop:partitions}
Let $m\ge 2$ with $m\not\equiv 2,3\pmod 5$. For $0\le t\le m$ we have
\[
[q^t]P^{(m)}(q)=\Delta^2p(t).
\]
Moreover, $\Delta^2 p(0)=1$, $\Delta^2 p(1)=-1$, and $\Delta^2 p(t)\ge 0$ for all $t\ge 2$. In particular, the coefficients of $q^0,\dots,q^{m}$ in $P^{(m)}$ are non-negative except for the coefficient $-1$ of $q^1$.
\end{proposition}

\begin{proof}
The coefficient of $q^t$ in $\qbin{6m}{m}$ is the number of partitions of $t$ with at most $m$ parts, each part at most $5m$. For $t\le m$ these restrictions are void, so $[q^t]G_m=p(t)$. Since $a=2m-1>m$, we have
\[
\frac{(1-q)^2}{(1-q^{a})(1-q^{b})}\equiv (1-q)^2 \pmod{q^{m+1}},
\]
and the first claim follows.

For the second claim, note that $p(t)-p(t-1)=p_1'(t)$ is the number of partitions of $t$ with no part equal to $1$ (obtained by removing one part $1$). Hence $\Delta^2p(t)=p'_1(t)-p'_1(t-1)$. We have $p'_1(0)=1$, $p'_1(1)=0$, $p'_1(2)=1$, which gives the values at $t=0,1,2$. For $t\ge 3$, adding $1$ to the largest part maps partitions of $t-1$ without parts $1$ injectively into partitions of $t$ without parts $1$. Hence $p'_1(t)\ge p'_1(t-1)$.
\end{proof}

Thus, the coefficient $-1$ of $q^1$ has a transparent origin: in the injection above, the empty partition of $0$ has no image. By reciprocity, the coefficient of $q^{D_m-1}$ is also $-1$. Conjecture~\ref{conj:GK} can therefore be read as a boxed analogue of the inequality $\Delta^2p(t)\ge 0$ ($t\ge2$), in which the box is $m\times 5m$ and the factor $(1-q)^2$ is replaced by $(1-q)^2/((1-q^{2m-1})(1-q^{3m-1}))$.

\subsection{The unboxed conjecture and the first $5m$ coefficients}\label{sec:tail}

We now extend Proposition~\ref{prop:partitions} from $t\le m$ to $t\le 5m$. For $t\le 5m$ the bound $5m$ on the size of the parts is void, but the bound $m$ on the number of parts is not. Instead of $p(t)$ we must therefore use partitions with at most $m$ parts. By conjugation, these are the partitions with parts at most $m$. The relevant inequality is the following.

For $k\ge2$, let $h_k(n)$ be the number of partitions of $n$ into parts from $\{2,3,\dots,k\}$ (with $h_k(n)=0$ for $n<0$), and put $d_k(n)=h_k(n)-h_k(n-1)$. Then
\begin{equation}\label{eq:dk}
\sum_{n\ge0}d_k(n)q^n=\frac{1-q}{(1-q^2)(1-q^3)\cdots(1-q^k)}=\frac{(1-q)^2}{(q;q)_k},
\end{equation}
where $(q;q)_k=(1-q)(1-q^2)\cdots(1-q^k)$. For $k\ge 7$ the first values are
\[
(d_k(0),\dots,d_k(7))=(1,-1,1,0,1,0,2,0).
\]

\begin{lemma}\label{lem:monotone}
Let $k\ge 9$. Then $h_k(n)\ge h_k(n-1)$ for every $n\ge2$, and $h_k(n)>h_k(n-1)$ for every $n\ge 8$. In other words, $d_k(n)\ge0$ for $n\ge 2$ and $d_k(n)\ge1$ for $n\ge 8$.
\end{lemma}

\begin{proof}
Let $\mathcal P(n)$ denote the set of partitions of $n$ into parts from $\{2,\dots,k\}$. We define a map $\iota\colon\mathcal P(n-1)\to\mathcal P(n)$. Let $\lambda\in\mathcal P(n-1)$.
\begin{enumerate}
\item[(1)] If $\lambda$ has a part $2$, replace one part $2$ by $3$.
\end{enumerate}
Otherwise, let $r$ be the number of parts $3$ in $\lambda$, and let $\nu$ be $\lambda$ with all parts $3$ removed; so $\nu$ has all its parts in $\{4,\dots,k\}$.
\begin{enumerate}
\item[(2a)] If $r$ is odd, let $\iota(\lambda)=\nu\cup 2^{(3r+1)/2}$.
\item[(2b)] If $r\ge2$ is even, let $\iota(\lambda)=\nu\cup 5\cup 2^{(3r-4)/2}$.
\item[(2c)] If $r=0$ and the largest part $\lambda_1$ satisfies $\lambda_1<k$, increase $\lambda_1$ by $1$.
\item[(2d)] If $r=0$ and $\lambda_1=k$, replace one part $k$ by the parts $k-5,2,2,2$.
\end{enumerate}
Here $\nu\cup 2^{j}$ means that $j$ parts equal to $2$ are added. In each case the parts of $\iota(\lambda)$ lie in $\{2,\dots,k\}$ (we use $5\le k$ and $k-5\ge4$), and $|\iota(\lambda)|=|\lambda|+1$, since $3r+1=2\cdot\frac{3r+1}2=5+2\cdot\frac{3r-4}{2}$ and $k+1=(k-5)+6$.

Let $\mu=\iota(\lambda)$, and let $j$ be the number of parts $2$ of $\mu$. In case (1), $\mu$ has a part $3$; in all other cases it has none. In case (2a), $j=(3r+1)/2\equiv2\pmod 3$; in case (2b), $j=(3r-4)/2\equiv1\pmod 3$; in case (2c), $j=0$; and in case (2d), $j=3$. Hence the five cases have pairwise disjoint images. Each case is injective: in case (1) we recover $\lambda$ by replacing a part $3$ of $\mu$ by $2$; in cases (2a) and (2b), $r$ is determined by $j$, and $\nu$ is obtained by deleting all parts $2$ (and, in case (2b), one part $5$); in case (2c), the largest part of $\mu$ is strictly larger than all other parts and we decrease it by $1$; in case (2d), we delete three parts $2$ and one part $k-5$ and insert a part $k$. Thus $\iota$ is injective and $h_k(n)\ge h_k(n-1)$.

Now let $n\ge 8$. Every integer $N\ge 6$ is a sum of elements of $\{4,6,7,8,9\}$: this holds for $6\le N\le 9$, and for $N\ge 10$ it follows by induction using $N=(N-4)+4$. Let $\rho$ be such a partition of $n-2$, and let $\mu=\rho\cup 2$. Then $\mu\in\mathcal P(n)$, $\mu$ has no part $3$, exactly one part $2$ (so $j=1$) and no part $5$. By the above description of the images, $\mu$ is not in the image of $\iota$ (it could only come from case (2b), whose images contain a part $5$). Hence $h_k(n)>h_k(n-1)$.
\end{proof}
A direct computation shows that $h_k(n)\ge h_k(n-1)$ for
$k\in\{5,7,8\}$ and $2\le n\le400$. Within this range, the strict
inequality $h_k(n)>h_k(n-1)$ for every $n\ge8$ holds for $k=7$,
but fails for $k=5$ at $n=9,11,13$ and for $k=8$ at $n=9$.

\begin{theorem}[the unboxed conjecture]\label{thm:unboxed}
Let $m\ge 9$, and let $a\ne b$ be integers with $a,b\ge 7$. Then the power series
\[
E_{a,b,m}(q)=\frac{(1-q)^2}{(1-q^a)(1-q^b)\,(q;q)_m}
\]
has constant term $1$, coefficient $-1$ at $q^1$, and non-negative coefficients at $q^t$ for all $t\ge 2$.
\end{theorem}

\begin{proof}
Write $1/((1-q^a)(1-q^b))=\sum_{n\ge0}N(n)q^n$, where $N(n)$ is the number of pairs $(i,j)\in\N^2$ with $ia+jb=n$. By \eqref{eq:dk} with $k=m$,
\[
[q^t]E_{a,b,m}=\sum_{n=0}^{t}N(n)\,d_m(t-n).
\]
For $t\le1$ this gives $1$ and $-1$. Let $t\ge2$. By Lemma~\ref{lem:monotone}, the only possibly negative term is the one with $t-n=1$, which equals $-N(t-1)$. Since $a+1\ge8$ and $b+1\ge8$, the terms with $n=t-1-a$ and $n=t-1-b$ (two distinct indices) are at least $N(t-1-a)$ and $N(t-1-b)$, respectively. Every representation $t-1=ia+jb$ has $(i,j)\ne(0,0)$, because $t-1\ge1$. If $i\ge1$ it yields the representation $(i-1,j)$ of $t-1-a$, and otherwise $j\ge1$ and it yields the representation $(i,j-1)$ of $t-1-b$. This map is injective, so $N(t-1)\le N(t-1-a)+N(t-1-b)$. Hence $[q^t]E_{a,b,m}\ge0$.
\end{proof}

Theorem~\ref{thm:unboxed} is the limiting case $N\to\infty$ of the Guo--Krattenthaler phenomenon. Indeed, $\qbin{N}{m}\to 1/(q;q)_m$ coefficientwise as $N\to\infty$, so
\[
\frac{(1-q)^2}{(1-q^{a})(1-q^{b})}\qbin{N}{m}\longrightarrow E_{a,b,m}(q).
\]
In our situation the box has width $5m$, and the limit is already reached in the first $5m+1$ coefficients.

\begin{corollary}\label{cor:tail}
Let $m\ge 9$ with $m\not\equiv2,3\pmod 5$. Then $[q^t]P^{(m)}=[q^t]E_{2m-1,3m-1,m}$, and this coefficient is non-negative, for $2\le t\le 5m$. By reciprocity, $[q^t]P^{(m)}\ge0$ also for $D_m-5m\le t\le D_m-2$.
\end{corollary}

\begin{proof}
We have $\qbin{6m}{m}=\prod_{j=1}^{m}(1-q^{5m+j})/(q;q)_m\equiv 1/(q;q)_m\pmod{q^{5m+1}}$. Multiplying by $(1-q)^2/((1-q^{2m-1})(1-q^{3m-1}))$ gives $P^{(m)}\equiv E_{2m-1,3m-1,m}\pmod{q^{5m+1}}$. Now apply Theorem~\ref{thm:unboxed} with $a=2m-1\ge17$ and $b=3m-1$.
\end{proof}

\begin{remark}\label{rem:exceptions}
Corollary~\ref{cor:tail} also explains the two exceptional values $m=4$ and $m=6$ in Theorem~\ref{thm:computation}. For $m=4$ we have $a=7$, and since $5<a$, $[q^5]P^{(4)}=d_4(5)=h_4(5)-h_4(4)=-1$. For $m=6$ we have $a=11$, and $[q^7]P^{(6)}=d_6(7)=-1$. In both cases the extra negative coefficient is precisely a failure of the monotonicity in Lemma~\ref{lem:monotone}. On the other hand, $d_6(13)=-1$ is compensated by $N(11)\,d_6(2)=1$, so $[q^{13}]P^{(6)}=0$.
\end{remark}

\begin{remark}
The naive approach would require the coefficient sequence of $\qbin{6m}{m}$ to be convex on the whole range $[2,D_m/2]$. This is false. For example, for $m=10$ the second difference first becomes negative at $t=199$, while $D_{10}/2=227$; asymptotically the inflection point of the Gaussian coefficients lies at distance of order $m^{3/2}$ from the centre, whereas $D_m/2$ is at distance of order $m$. The reduction below avoids convexity and uses only a monotonicity property after one division.
\end{remark}

\section{Andrews' lemma with one defect, and a reduction}\label{sec:reduction}

The following lemma is an analogue of \cite[Lemma 5.1]{GK}, in which the unimodal input polynomial is allowed to have one ``dip'' at $q^1$.

\begin{lemma}\label{lem:defect}
Let $U(q)=\sum_{s}u_sq^s\in\Z[q]$ be reciprocal, let $b\ge 2$, and suppose that $P(q)=U(q)/[b]_q$ is a polynomial of degree $D\ge 3$. Assume that
\begin{enumerate}
\item[(i)] $u_0=1$ and $u_1=0$;
\item[(ii)] $u_{s}\ge u_{s-1}$ for $2\le s\le \lfloor D/2\rfloor$;
\item[(iii)] if $b+1\le \lfloor D/2\rfloor$, then $u_{b+1}>u_b$.
\end{enumerate}
Then $[q^1]P=[q^{D-1}]P=-1$, and all other coefficients of $P$ are non-negative.
\end{lemma}

\begin{proof}
Put $\delta_s=u_s-u_{s-1}$ (with $u_{-1}=0$), so that $\delta_0=1$, $\delta_1=-1$, and, by (ii), $\delta_s\ge 0$ for $2\le s\le\lfloor D/2\rfloor$. From $P(q)(1-q^b)=U(q)(1-q)$ we get $p_t-p_{t-b}=\delta_t$, hence
\[
p_t=\sum_{j\ge 0,\ jb\le t}\delta_{t-jb}\qquad (t\ge 0).
\]
Let $0\le t\le\lfloor D/2\rfloor$. All indices $t-jb$ lie in $[0,\lfloor D/2\rfloor]$. If $t\not\equiv 1\pmod b$, then no index equals $1$ and $p_t\ge 0$. If $t=1$, then $p_1=\delta_1=-1$ because $b\ge 2$. If $t=1+jb$ with $j\ge 1$, the sum contains $\delta_1=-1$ together with $\delta_{b+1}$, and $b+1\le t\le\lfloor D/2\rfloor$, so $\delta_{b+1}\ge 1$ by (iii). All other terms are non-negative, hence $p_t\ge 0$.

Since $U$ and $[b]_q$ are reciprocal, so is $P$, i.e.\ $p_{D-t}=p_t$. As $D\ge 3$, we have $D-1>\lfloor D/2\rfloor\ge 1$, and the statement follows.
\end{proof}

Note that the hypothesis is only a \emph{one-sided} monotonicity on the first half of $P$; it does not require $U$ to be unimodal, and it tolerates the dip $1,0$ at the beginning.

\begin{theorem}\label{thm:reduction}
Let $m\ge 2$ with $m\not\equiv 2,3\pmod 5$, and let
\[
U_m(q)=\frac{1}{[2m-1]_q}\qbin{6m}{m},\qquad V_m(q)=\frac{1}{[3m-1]_q}\qbin{6m}{m}.
\]
Suppose that one of the following holds:
\begin{enumerate}
\item[(H$_U$)] the coefficients $u_s$ of $U_m$ satisfy $u_s\ge u_{s-1}$ for $2\le s\le \lfloor D_m/2\rfloor$ and $u_{3m}>u_{3m-1}$;
\item[(H$_V$)] the coefficients $v_s$ of $V_m$ satisfy $v_s\ge v_{s-1}$ for $2\le s\le \lfloor D_m/2\rfloor$ and $v_{2m}>v_{2m-1}$.
\end{enumerate}
Then $[q^1]P^{(m)}=[q^{D_m-1}]P^{(m)}=-1$, and all other coefficients of $P^{(m)}$ are non-negative.
\end{theorem}

\begin{proof}
By Proposition~\ref{prop:poly}, $P^{(m)}$ is a polynomial, hence so are $U_m=[3m-1]_qP^{(m)}$ and $V_m=[2m-1]_qP^{(m)}$. Both are reciprocal. Since $[q^0]G_m=[q^1]G_m=1$, we have $u_0=v_0=1$ and $u_1=v_1=0$. Apply Lemma~\ref{lem:defect} to $U=U_m$, $b=3m-1$ (respectively $U=V_m$, $b=2m-1$). Note that $D_m\ge 3$ for $m\ge 2$.
\end{proof}

In particular, Conjecture~\ref{conj:GK} for a given $n$ follows from (H$_U$) or (H$_V$) for $m=5n$. The conditions (H$_U$) and (H$_V$) concern a quotient by a \emph{single} $q$-integer, which is exactly the situation of Andrews' lemma.

\begin{theorem}\label{thm:computation}
For every $m$ with $2\le m\le 300$ and $m\not\equiv 2,3\pmod 5$:
\begin{enumerate}
\item if $m\notin\{4,6\}$, then (H$_U$) holds; moreover $U_m$ is weakly increasing from $q^1$ up to its middle coefficient;
\item if $m\notin\{4,5,6\}$, then (H$_V$) also holds;
\item for $m\in\{4,6\}$, the polynomial $P^{(m)}$ has four negative coefficients (each equal to $-1$), namely at $q^1,q^5,q^{59},q^{63}$ for $m=4$ and at $q^1,q^7,q^{147},q^{153}$ for $m=6$.
\end{enumerate}
Consequently, Conjecture~\ref{conj:GK} holds for all $n\le 60$, and $P^{(m)}$ has exactly the two exceptional coefficients $-1$ for all $m\le 300$ with $m\not\equiv 2,3\pmod 5$ and $m\notin\{1,4,6\}$.
\end{theorem}

\begin{proof}
This is a finite computation with exact integer arithmetic. The polynomial $\qbin{N}{k}$ is built as $\prod_{i=0}^{k-1}(1-q^{N-i})/(1-q^{i+1})$, alternately multiplying by $1-q^{r}$ and dividing exactly by $1-q^{s}$; each partial product is the Gaussian polynomial $\qbin{N}{i+1}$, so all divisions are exact. Division by $[c]_q$ is performed as multiplication by $1-q$ followed by exact division by $1-q^c$, using the recurrence $x_t=y_t+x_{t-c}$, and we check that the remainder vanishes. For each $m$ we computed $G_m$, $U_m$, $V_m$ and $P^{(m)}$, checked (H$_U$) and (H$_V$) directly, and listed the negative coefficients of $P^{(m)}$. The largest case $m=300$ involves polynomials of degree about $4.5\cdot 10^5$. For $m=1$ we have $P^{(1)}=1+q^2+q^4$ (here $2m-1=1$).
\end{proof}

The data suggest the following extension of Conjecture~\ref{conj:GK}.

\begin{conjecture}\label{conj:general}
Let $m\ge 9$ (or $m=5$) with $m\not\equiv 2,3\pmod 5$. Then the polynomial
\[
\frac{(1-q)^2}{(1-q^{2m-1})(1-q^{3m-1})}\qbin{6m}{m}
\]
has non-negative coefficients, except for the coefficients of $q^1$ and $q^{5m^2-5m+3}$, which are equal to $-1$.
\end{conjecture}

\begin{conjecture}\label{conj:almost}
For every $m\ge 9$ (or $m=5$) with $m\not\equiv 2,3\pmod 5$, the polynomial $\qbin{6m}{m}/[2m-1]_q$ is weakly increasing from its coefficient of $q^1$ up to its middle coefficient.
\end{conjecture}

By Theorem~\ref{thm:reduction}, Conjecture~\ref{conj:almost} (together with the easy inequality $u_{3m}>u_{3m-1}$) implies Conjectures~\ref{conj:general} and~\ref{conj:GK}. Conjecture~\ref{conj:almost} has the same flavour as the unimodality of Gaussian polynomials; strict unimodality results such as those of Pak and Panova \cite{PP} indicate that it might be approachable by similar methods. In Section~\ref{sec:asymptotic} we bypass Conjecture~\ref{conj:almost} and prove Conjecture~\ref{conj:general} directly for all $m\ge10^{30}$.

\section{An explicit bound: the conjecture holds for $n\ge 2\cdot10^{29}$}\label{sec:asymptotic}

In this section we prove that Conjecture~\ref{conj:general}, and hence Conjecture~\ref{conj:GK}, holds for all parameters beyond an explicit threshold. Unlike Section~\ref{sec:reduction}, we analyse $P^{(m)}$ directly and do not go through $U_m$. The reason is that the coefficients of $P^{(m)}$ involve no cancellation at leading order. The factor $(1-q)^2/((1-q^{a})(1-q^{b}))$ is smooth near the relevant saddle point, so $[q^t]P^{(m)}$ is a \emph{positive} multiple of a local central limit probability, up to a controlled error. For $U_m$, by contrast, monotonicity is a second-order effect.

\begin{theorem}\label{thm:asymptotic}
Let $m\ge 10^{30}$ with $m\not\equiv 2,3\pmod 5$. Then all coefficients of $P^{(m)}(q)$ are non-negative, except $[q^1]P^{(m)}=[q^{D_m-1}]P^{(m)}=-1$. In particular, Conjecture~\ref{conj:GK} holds for all $n\ge 2\cdot 10^{29}$.
\end{theorem}

No attempt has been made to optimise the threshold; see Remark~\ref{rem:numerics}. All constants below are explicit, and $\log$ denotes the natural logarithm. We fix $m\ge10^{30}$ with $m\not\equiv2,3\pmod5$, and write $N=6m$, $a=2m-1$, $b=3m-1$, $G=G_m$, $P=P^{(m)}=\sum_tp_tq^t$, $D=D_m$ and $T=\lfloor D/2\rfloor$. By Proposition~\ref{prop:partitions} and Corollary~\ref{cor:tail}, $p_1=-1$ and $p_t\ge0$ for $2\le t\le 5m$. Since $P$ is reciprocal, and $t\mapsto D-t$ maps $(T,D-5m)$ into $(5m,T]$, Theorem~\ref{thm:asymptotic} follows from
\begin{equation}\label{eq:goal}
p_t>0\qquad\text{for all } 5m<t\le T .
\end{equation}
From now on we fix such a $t$.

\subsection{A probabilistic representation}
Let
\[
F(w,z)=\prod_{j=0}^{N-1}(1+wz^j)=\sum_{S\subseteq\{0,\dots,N-1\}}w^{|S|}z^{\sum_{s\in S}s},
\qquad
K(z)=\frac{(1-z)^2}{(1-z^a)(1-z^b)}=\frac{1}{[a]_z[b]_z} .
\]
Then $[w^m]F(w,z)=z^{\binom m2}G(z)$. The function $K$ is analytic and zero-free in $|z|<1$, and $P=GK$ as power series. Hence, with $t'=t+\binom m2$, for all $\rho>0$ and $0<r<1$,
\begin{equation}\label{eq:cauchy}
p_t=[w^mz^{t'}]F(w,z)K(z)=\frac{1}{4\pi^2}\int_{-\pi}^{\pi}\!\!\int_{-\pi}^{\pi}F(\rho e^{i\psi},re^{i\varphi})K(re^{i\varphi})\rho^{-m}r^{-t'}e^{-i(m\psi+t'\varphi)}\,d\psi\,d\varphi .
\end{equation}
For $\lambda\in\mathbb R$ and $\theta>0$ put $\rho=e^{\lambda}$, $r=e^{-\theta}$, $\mu=\lambda/\theta$, $s=1/\theta$, $A=1/(1-e^{-\theta})$, and
\[
p_j=\frac{1}{1+e^{\theta(j-\mu)}},\qquad v_j=p_j(1-p_j),\qquad q_j=\min(p_j,1-p_j)\qquad(0\le j<N).
\]
Note that $s\le A\le s+1$, and that $v_j\le q_j\le e^{-\theta|j-\mu|}$. Then
\[
\frac{F(\rho e^{i\psi},re^{i\varphi})}{F(\rho,r)}=\prod_{j=0}^{N-1}\bigl(1-p_j+p_je^{i(\psi+j\varphi)}\bigr)=:\Phi(\psi,\varphi)
\]
is the characteristic function of $(\sum_jB_j,\sum_j jB_j)$, where the $B_j$ are independent Bernoulli variables with $\Pr(B_j=1)=p_j$. We write
\[
V=\sum_j v_j,\qquad \bar\jmath=\frac1V\sum_j jv_j,\qquad W=\sum_j v_j(j-\bar\jmath)^2,\qquad M_3=\sum_j v_j|j-\bar\jmath|^3 .
\]
Since $\bar\jmath$ minimises $x\mapsto\sum_jv_j(j-x)^2$ and $v_j\le\frac14$, we always have
\begin{equation}\label{eq:trivial}
V\le\tfrac N4=1.5m,\qquad W\le\tfrac14\sum_{j=0}^{N-1}\Bigl(j-\tfrac{N-1}2\Bigr)^2\le\tfrac{N^3}{48}=4.5m^3,\qquad M_3\le\tfrac14\sum_{j=0}^{N-1}j^3\le\tfrac{N^4}{16}=81m^4 .
\end{equation}
(For $M_3$ we use that $x\mapsto\sum_{j<N}|j-x|^3$ is convex, so its maximum on $[0,N-1]$ is attained at an endpoint.)

\begin{lemma}\label{lem:mu}
Suppose that $\sum_jp_j=m$ and $A<m-1$. Then $|\mu-m|\le A+1$.
\end{lemma}

\begin{proof}
Using $p_j\le e^{-\theta(j-\mu)}$ for $j>\mu$, we get $m\le(\lfloor\mu\rfloor+1)^++A$. Since $A<m-1$, this forces $\mu\ge0$ and $\mu\ge m-1-A$. Using $1-p_j\le e^{-\theta(\mu-j)}$ for $j<\mu$, we get $m\ge\min(\mu,N)-A$. Since $A<5m$, this forces $\mu<N$ and $\mu\le m+A$.
\end{proof}

\begin{lemma}[saddle point]\label{lem:saddle}
There is a unique pair $(\lambda,\theta)$ with $\theta>0$ such that $\sum_jp_j=m$ and $\sum_jjp_j=t'$. It satisfies
\[
\frac{4}{9m^2}\le\theta<\frac{1.06}{\sqrt m},\qquad\text{so that}\qquad 0.9434\sqrt m<s\le 2.25\,m^2 .
\]
\end{lemma}

\begin{proof}
For fixed $\theta\ge0$, the map $\lambda\mapsto\sum_jp_j$ is an increasing bijection $\mathbb R\to(0,N)$, so $\sum_jp_j=m$ defines a smooth function $\lambda(\theta)$. Let $E(\theta)=\sum_jjp_j$ at $\lambda=\lambda(\theta)$. Since $\partial p_j/\partial\lambda=v_j$ and $\partial p_j/\partial\theta=-jv_j$ (as functions of $(\lambda,\theta)$), differentiating the constraint gives $\lambda'(\theta)=\bar\jmath$, and then
\[
E'(\theta)=\sum_j jv_j(\bar\jmath-j)=-W<0 .
\]
At $\theta=0$ all $p_j$ equal $1/6$ and $E(0)=\binom m2+\frac52m^2$; as $\theta\to\infty$, $E(\theta)\to\sum_{j<m}j=\binom m2$. Since $\binom m2<t'<E(0)$, there is a unique solution.

Lower bound. By \eqref{eq:trivial}, $W\le4.5m^3$ for every $\theta$, so $\frac52m^2-t=E(0)-E(\theta)\le4.5m^3\theta$. Since $t\le T\le\frac12(5m^2-5m+4)$, the left-hand side is at least $2m$, and $\theta\ge 4/(9m^2)$.

Upper bound. Let $\tau(\theta)=E(\theta)-\binom m2$, which is decreasing. Using $\sum_jp_j=m$, one checks that
\[
\tau(\theta)=\sum_{j\ge m}(j-m)p_j+\sum_{j<m}(m-j)(1-p_j).
\]
Let $\theta_1=1.06/\sqrt m$ and let $A_1$ be the corresponding value of $A$. Then $A_1\le 1/\theta_1+1<m-1$, so by Lemma~\ref{lem:mu}, $|\mu-m|\le A_1+1$ at $\theta=\theta_1$. Hence $p_j\le e^{\theta_1(A_1+1)}e^{-\theta_1(j-m)}$ for $j\ge m$, and $1-p_j\le e^{\theta_1(A_1+1)}e^{-\theta_1(m-j)}$ for $j<m$. Using $\sum_{i\ge0}ie^{-\theta i}=e^{-\theta}A^2$, we obtain
\[
\tau(\theta_1)\le 2e^{\theta_1A_1}A_1^2\le\frac{2e^{1+\theta_1}(1+\theta_1)^2}{\theta_1^2}\le 4.86\,m<5m<t ,
\]
because $\theta A\le1+\theta$, $A\le(1+\theta)/\theta$ and $\theta_1\le 0.0011$. Since $\tau$ is decreasing and $\tau(\theta)=t$, we get $\theta<\theta_1$.
\end{proof}

From now on, $(\lambda,\theta)$ is the saddle point of Lemma~\ref{lem:saddle}. We distinguish two regimes:
\[
\text{(I)}\ \ \theta\ge 22/m\ \ (\text{so } 0.9434\sqrt m<s\le m/22),\qquad\qquad \text{(II)}\ \ \theta<22/m\ \ (\text{so } s>m/22).
\]
Put $\Lambda=\log 11$ and $J=\{j\in[0,N):\,|j-\mu|\le\Lambda s\}=\{j:\frac1{12}\le p_j\le\frac{11}{12}\}$. Then $v_j\ge\frac{11}{144}$ on $J$.

\begin{lemma}\label{lem:sums}
For every real $\mu$, every $\theta>0$ and $k=0,1,2,3$,
\[
\sum_{j\in\mathbb Z}|j-\mu|^ke^{-\theta|j-\mu|}\le c_kA^{k+1},\qquad (c_0,c_1,c_2,c_3)=(2,2,4,12).
\]
\end{lemma}

\begin{proof}
On each side of $\mu$ the distances are $u,u+1,u+2,\dots$ with $0\le u<1$ (respectively $0<u\le1$), so the sum is at most $2\sum_{n\ge0}(n+1)^ke^{-\theta n}=2e^{\theta}\sum_{n\ge1}n^ke^{-\theta n}$. Now use $\sum_{n\ge1}x^n=xA$, $\sum nx^n=xA^2$, $\sum n^2x^n=x(1+x)A^3$ and $\sum n^3x^n=x(1+4x+x^2)A^4$, with $x=e^{-\theta}$.
\end{proof}

\begin{lemma}[shape of the tilted measure]\label{lem:shape}
\begin{enumerate}
\item[(I)] In regime (I), $J$ consists of consecutive integers, $|J|\ge 4.75s$, and
\[
0.3628\,s\le V\le 2.003\,s,\quad 0.682\,s^3\le W\le 4.02\,s^3,\quad |\bar\jmath-\mu|\le 3.33\,s,\quad M_3\le 193\,s^4 .
\]
\item[(II)] In regime (II), $J$ consists of consecutive integers, $|J|\ge 0.217m$, and
\[
0.01658\,m\le V\le 1.5\,m,\qquad 6.5\cdot10^{-5}m^3\le W\le4.5\,m^3,\qquad M_3\le 81\,m^4 .
\]
\end{enumerate}
\end{lemma}

\begin{proof}
(I) Since $A\le s+1<m-1$, Lemma~\ref{lem:mu} gives $|\mu-m|\le s+2$. Hence $[\mu-\Lambda s,\mu+\Lambda s]\subseteq[m-3.4s-2,\,m+3.4s+2]\subseteq[0.845m-2,\,1.155m+2]\subseteq[0,N-1]$, so $J$ is the set of all integers in this window and $|J|\ge 2\Lambda s-1\ge4.75s$. Thus $V\ge\frac{11}{144}|J|\ge0.3628s$. By Lemma~\ref{lem:sums} and $v_j\le e^{-\theta|j-\mu|}$, we have $V\le2A$, $S_1:=\sum_jv_j|j-\mu|\le2A^2$, $S_2:=\sum_jv_j(j-\mu)^2\le4A^3$ and $S_3:=\sum_jv_j|j-\mu|^3\le12A^4$. Since $A\le s+1\le1.00107s$ (as $s\ge 943$), this gives $V\le2.003s$, $S_1\le2.005s^2$, $S_2\le4.014s^3$ and $S_3\le12.06s^4$. Then $W\le S_2\le4.02s^3$. By Cauchy--Schwarz, $|\bar\jmath-\mu|^2\le S_2/V$, so $d:=|\bar\jmath-\mu|\le3.33s$. Since $J$ consists of $|J|$ consecutive integers, $W\ge\frac{11}{144}\cdot\frac{|J|^3-|J|}{12}\ge0.682s^3$. Finally,
\[
M_3\le\sum_jv_j(|j-\mu|+d)^3=S_3+3dS_2+3d^2S_1+d^3V\le193s^4 .
\]

(II) Since $p_j$ is decreasing in $j$, $J$ consists of consecutive integers. It is non-empty: not all $p_j$ are $<\frac1{12}$ (otherwise $\sum p_j<m/2$), not all are $>\frac{11}{12}$ (otherwise $\sum p_j>m$), and the logit $\theta(\mu-j)$ decreases by $\theta<1<2\Lambda$ from one $j$ to the next. If $J$ contains neither $0$ nor $N-1$, then $|J|\ge2\Lambda s-1\ge 4.79m/22-1$. If $0\in J$, then $p_j<\frac1{12}$ for $j\ge|J|$, so $m\le|J|+\frac{N}{12}$ and $|J|\ge m/2$. If $N-1\in J$, then $p_j>\frac{11}{12}$ for $j<N-|J|$, so $m\ge\frac{11}{12}(N-|J|)$ and $|J|\ge 4.9m$. In all cases $|J|\ge0.217m$. Hence $V\ge\frac{11}{144}\cdot0.217m\ge0.01658m$ and $W\ge\frac{11}{144}\cdot\frac{|J|^3-|J|}{12}\ge6.5\cdot10^{-5}m^3$. The upper bounds are \eqref{eq:trivial}.
\end{proof}

\begin{lemma}[the factor $K$]\label{lem:K}
Let $z=re^{i\varphi}$.
\begin{enumerate}
\item[(i)] $K(r)>0$ and $|K(z)|\le 4A^2K(r)$.
\item[(ii)] $\beta:=rK'(r)/K(r)$ satisfies $-\min(2s,5m)\le\beta\le0$.
\item[(iii)] Let $f(\varphi)=\log K(re^{i\varphi})$, continued along the arc from $\varphi=0$. Then $f'(0)=i\beta$, and $|f(\varphi)-f(0)-i\beta\varphi|\le\frac12H\varphi^2$, where $H=4s^2$ in general, and $H=4.76m^2$ if $\theta<\frac1{6m}$ and $|\varphi|\le\frac1{6m}$.
\end{enumerate}
\end{lemma}

\begin{proof}
(i) We have $|1-z|\le2$ and $|1-z^n|\ge1-r^n$, so $|K(z)|\le4/((1-r^a)(1-r^b))=4A^2K(r)$.

(ii) We have $\beta=-\frac{2r}{1-r}+\frac{ar^a}{1-r^a}+\frac{br^b}{1-r^b}$, and $0\le\frac{nr^n}{1-r^n}\le\frac{r}{1-r}\le s$. Also $\beta=-r\frac{d}{dr}\bigl(\log[a]_r+\log[b]_r\bigr)\ge-(a-1)-(b-1)\ge-5m$.

(iii) We have $f''(\varphi)=-\bigl((z\tfrac{d}{dz})^2\log K\bigr)(z)$ and $(z\frac d{dz})^2\log(1-z^n)=-n^2z^n/(1-z^n)^2$. Hence, for $|z|=r$,
\[
|f''|\le\frac{2r}{(1-r)^2}+\sum_{n\in\{a,b\}}\frac{n^2r^n}{(1-r^n)^2}=\frac{2}{4\sinh^2(\theta/2)}+\sum_{n\in\{a,b\}}\frac{n^2}{4\sinh^2(n\theta/2)}\le\frac{4}{\theta^2}=4s^2 .
\]
Now let $\theta<\frac1{6m}$ and $|\varphi|\le\frac1{6m}$. Then $|1-z|\le\theta+|\varphi|<\frac1b$. Write $\log[n]_z=\sum_{k=1}^{n-1}\log(1-z/\zeta_k)$ with $\zeta_k=e^{2\pi ik/n}$ and $n\in\{a,b\}$. For $1\le k\le n/2$ we have $|\zeta_k-1|=2\sin(\pi k/n)\ge4k/n$, hence $|\zeta_k-z|\ge3k/n$, and symmetrically for $n-k$. Since $(z\frac d{dz})^2\log(1-z/\zeta)=-z\zeta/(\zeta-z)^2$, we get
\[
|f''|\le\sum_{n\in\{a,b\}}2\sum_{k\ge1}\frac{n^2}{9k^2}\le\frac{\pi^2}{27}(a^2+b^2)\le4.76m^2 .
\]
Taylor's formula with integral remainder completes the proof.
\end{proof}

\subsection{Major and minor arcs}
We change variables from $\psi$ to $\alpha\in[-\pi,\pi)$, $\alpha\equiv\psi+\bar\jmath\varphi\pmod{2\pi}$. This does not change the integral in \eqref{eq:cauchy}. Put $\tilde s_j=\alpha+(j-\bar\jmath)\varphi$. Then $\tilde s_j\equiv\psi+j\varphi\pmod{2\pi}$, with the same multiple of $2\pi$ for all $j$. Since $m$ is an integer, this gives
\[
e^{-i(m\psi+t'\varphi)}=e^{-i\sum_jp_j\tilde s_j},\qquad \sum_jv_j\tilde s_j^2=V\alpha^2+W\varphi^2=:Q .
\]
Let
\[
L=\sqrt{24+13.8\log m},\qquad \kappa=LV^{-1/2},\qquad \eta=LW^{-1/2},\qquad \mathcal B=\{(\alpha,\varphi):|\alpha|\le\kappa,\ |\varphi|\le\eta\}.
\]
For $m\ge10^{30}$ we have $31.2\le L$ and $L\le 10^{-5}m^{1/4}$. Combined with Lemma~\ref{lem:shape}, this gives $\kappa\le\frac12$ and $\eta\le\frac1{30.66\,s}$ in regime (I), and $\eta\le\frac1{24m}$ in regime (II).

\begin{lemma}[major arc]\label{lem:major}
For $(\alpha,\varphi)\in\mathcal B$,
\[
\Phi(\psi,\varphi)\,e^{-i(m\psi+t'\varphi)}=\exp\bigl(-\tfrac12Q+R\bigr),\qquad |R|\le0.99\bigl(V|\alpha|^3+M_3|\varphi|^3\bigr)+R_{\mathrm{tail}},
\]
where $R_{\mathrm{tail}}=0$ in regime (II), and in regime (I)
\[
R_{\mathrm{tail}}\le 6m\cdot27.9\,\bigl(5+(\tfrac12+7.27Lms^{-3/2})^2\bigr)\exp\bigl(-0.4129\sqrt s/L\bigr).
\]
\end{lemma}

\begin{proof}
We have $\Phi e^{-i(m\psi+t'\varphi)}=\prod_j(1-p_j+p_je^{i\tilde s_j})e^{-ip_j\tilde s_j}$. For each $j$ we write the factor as $\exp(-\frac12v_j\tilde s_j^2+r_j)$.

If $|\tilde s_j|\le1$, put $g(s)=\log(1-p_j+p_je^{is})$. Then $g(0)=0$, $g'(0)=ip_j$, $g''(0)=-v_j$, and
\[
|g'''(s)|=\frac{v_j|1-p_j-p_je^{is}|}{|1-p_j+p_je^{is}|^3}\le\frac{v_j}{(1-\frac12(1-\cos1))^{3/2}}\le1.479\,v_j ,
\]
because $|1-p+pe^{is}|^2=1-2v(1-\cos s)$ and $v\le\frac14$. Hence $|r_j|\le0.2466\,v_j|\tilde s_j|^3$.

If $q_j\le\frac14$ and $s\in\mathbb R$ is arbitrary, write $1-p_j+p_je^{is}=1+w$ with $w=p_j(e^{is}-1)$ when $p_j\le\frac12$, and $=e^{is}(1+w)$ with $w=(1-p_j)(e^{-is}-1)$ when $p_j>\frac12$. In both cases $|w|\le2q_j\le\frac12$, so $|\log(1+w)|\le2|w|$, and we obtain $|r_j|\le q_j(4+|s|+\frac12s^2)\le q_j(5+s^2)$.

Let $\mathcal T=\{j:|j-\bar\jmath|>1/(2\eta)\}$. For $j\notin\mathcal T$ we have $|\tilde s_j|\le\kappa+\frac12\le1$. Using $|x+y|^3\le4(|x|^3+|y|^3)$,
\[
\sum_{j\notin\mathcal T}|r_j|\le0.2466\cdot4\sum_jv_j\bigl(|\alpha|^3+|j-\bar\jmath|^3|\varphi|^3\bigr)\le0.99\bigl(V|\alpha|^3+M_3|\varphi|^3\bigr).
\]
In regime (II), $1/(2\eta)=\sqrt W/(2L)\ge0.004m^{3/2}/L>N$, so $\mathcal T=\emptyset$. In regime (I), $j\in\mathcal T$ implies $|j-\bar\jmath|>0.4129\,s^{3/2}/L$, so that
\[
q_j\le e^{-\theta|j-\mu|}\le e^{\theta|\bar\jmath-\mu|}e^{-\theta|j-\bar\jmath|}\le27.9\,e^{-0.4129\sqrt s/L}\le\tfrac14 .
\]
Also $|\tilde s_j|\le\frac12+N\eta\le\frac12+7.27Lms^{-3/2}$. Summing $q_j(5+\tilde s_j^2)$ over at most $N=6m$ indices gives the bound for $R_{\mathrm{tail}}$.
\end{proof}

\begin{lemma}[minor arcs]\label{lem:minor}
For $(\alpha,\varphi)\in[-\pi,\pi)^2\setminus\mathcal B$ we have
\[
|\Phi(\psi,\varphi)|\le\exp\bigl(-\min(0.4367L^2,\,E_{\mathrm{far}})\bigr),
\]
where $E_{\mathrm{far}}=2.28\cdot10^{-6}s$ in regime (I) and $E_{\mathrm{far}}=3.5\cdot10^{-10}m$ in regime (II).
\end{lemma}

\begin{proof}
From $|1-p+pe^{is}|^2=1-2v(1-\cos s)$ we get $|\Phi|\le e^{-\mathcal E}$, where $\mathcal E=\sum_jv_j(1-\cos\tilde s_j)$. Outside $\mathcal B$ we have $Q>L^2$. Put $\varphi^*=\frac1{30.66\,s}$ in regime (I) and $\varphi^*=\frac1{24m}$ in regime (II).

\emph{Near region: $|\alpha|\le\frac12$ and $|\varphi|\le\varphi^*$.} Let $\Delta=15.33s$ in regime (I) and $\Delta=N$ in regime (II), and let $C=\{j:|j-\bar\jmath|\le\Delta\}$. For $j\in C$ we have $|\tilde s_j|\le\frac12+\Delta\varphi^*\le1$, and $1-\cos x\ge(1-\cos1)x^2\ge0.45969x^2,$ for $|x|\le1$. Hence
\[
\mathcal E\ge0.4597\Bigl(Q-\sum_{j\notin C}v_j\tilde s_j^2\Bigr)\ge0.4597\bigl(Q-2\alpha^2V_t-2\varphi^2W_t\bigr),
\]
where $V_t=\sum_{j\notin C}v_j$ and $W_t=\sum_{j\notin C}v_j(j-\bar\jmath)^2$. In regime (II), $C$ contains all $j$, so $V_t=W_t=0$. In regime (I), $j\notin C$ implies $|j-\mu|>\Delta-3.33s=12s$. Arguing as in Lemma~\ref{lem:sums}, we get $V_t\le2e^{-12}A\le V/40$ and
\[
W_t\le2e^{-12}\sum_{n\ge0}(15.34s+n)^2e^{-\theta n}\le2e^{-12}\bigl(15.34^2s^2A+30.68sA^2+2A^3\bigr)\le W/40 .
\]
In both regimes, therefore, $\mathcal E\ge0.4597\cdot0.95\,Q\ge0.4367L^2$.

\emph{Far region.} Let $J$ be as in Lemma~\ref{lem:shape}, put $\nu=\frac{11}{144}$, and recall $v_j\ge\nu$ on $J$.

(F1) If $|\varphi|\ge2\pi/|J|$, then $|\sum_{j\in J}e^{i(\psi+j\varphi)}|\le1/|\sin(\varphi/2)|\le\pi/|\varphi|\le|J|/2$. Hence $\mathcal E\ge\nu|J|/2$.

(F2) Let $\varphi^*<|\varphi|<2\pi/|J|$. Let $j_c$ be a middle element of $J$, let $h=\lfloor(|J|-1)/8\rfloor$, and let $\gamma\in[-\pi,\pi)$ with $\gamma\equiv\psi+j_c\varphi$. For $|k|\le h$ we have $j_c+k\in J$ and $|k\varphi|\le\pi/4$. If $|\gamma|\ge\pi/2$, then each $\gamma+k\varphi$ has distance at least $\pi/4$ from $2\pi\mathbb Z$, and $\mathcal E\ge\nu(1-\cos\frac\pi4)(2h+1)$. If $|\gamma|<\pi/2$, then $|\gamma+k\varphi|<3\pi/4$, and $1-\cos x\ge2x^2/\pi^2$ for $|x|\le\pi$ gives
\[
\mathcal E\ge\frac{2\nu}{\pi^2}\sum_{|k|\le h}(\gamma+k\varphi)^2\ge\frac{2\nu}{\pi^2}\cdot\frac{h(h+1)(2h+1)}{3}\,\varphi^{*2}.
\]

(F3) Let $|\varphi|\le\varphi^*$ and $|\alpha|>\frac12$. For $j\in J$ we have $|(j-\bar\jmath)\varphi|\le0.187$ in regime (I), using $|j-\bar\jmath|\le(\Lambda+3.33)s$, and $\le\frac14$ in regime (II). Hence $\tilde s_j$ has distance at least $\frac14$ from $2\pi\mathbb Z$, and $\mathcal E\ge\nu(1-\cos\frac14)|J|$.

Inserting $|J|\ge4.75s$ and $h\ge0.5925s$ (regime (I)), respectively $|J|\ge0.217m$ and $h\ge0.0271m$ (regime (II)), the smallest of these lower bounds is the second bound in (F2). It equals at least $2.28\cdot10^{-6}s$ in regime (I) and $3.5\cdot10^{-10}m$ in regime (II).
\end{proof}

\subsection{Proof of Theorem~\ref{thm:asymptotic}}
By \eqref{eq:cauchy},
\[
p_t=F(\rho,r)\rho^{-m}r^{-t'}K(r)\cdot I,\qquad I=\frac1{4\pi^2}\iint_{[-\pi,\pi)^2}\Phi(\psi,\varphi)e^{-i(m\psi+t'\varphi)}\frac{K(re^{i\varphi})}{K(r)}\,d\alpha\,d\varphi,
\]
and the prefactor is positive, so it suffices to show $I>0$. Split $I=I_{\mathcal B}+I_{\mathrm{minor}}$.

\emph{Major arc.} On $\mathcal B$, Lemmas~\ref{lem:major} and~\ref{lem:K}(iii) give the integrand $\exp(-\frac12Q+i\beta\varphi+X)$ with $X=R+f(\varphi)-f(0)-i\beta\varphi$. In regime (II) we use the second value of $H$ when $\theta<\frac1{6m}$ (note $\eta\le\frac1{24m}$), and $H=4s^2\le144m^2$ otherwise. This gives $|X|\le\varepsilon$, where
\begin{align*}
\text{(I)}\quad\varepsilon&=0.99L^3\Bigl(\frac{1}{\sqrt{0.3628}}+\frac{193}{0.682^{3/2}}\Bigr)s^{-1/2}+\frac{2L^2}{0.682\,s}+R_{\mathrm{tail}}\\
&\le 341.0\,L^3s^{-1/2}+2.94L^2s^{-1}+R_{\mathrm{tail}},\\
\text{(II)}\quad\varepsilon&=0.99L^3\Bigl(\frac{1}{\sqrt{0.01658}}+\frac{81}{(6.5\cdot10^{-5})^{3/2}}\Bigr)m^{-1/2}+\frac{72L^2}{6.5\cdot10^{-5}m}\\
&\le1.531\cdot10^{8}L^3m^{-1/2}+1.11\cdot10^6L^2m^{-1}.
\end{align*}
Let $\mathcal M=\frac1{4\pi^2}\iint_{\mathbb R^2}e^{-\frac12Q+i\beta\varphi}\,d\alpha\,d\varphi=\frac{e^{-\beta^2/(2W)}}{2\pi\sqrt{VW}}$. If $\varepsilon\le\frac12$, then $|e^X-1|\le1.65|X|$ and $\Pr(|Z|>L)\le e^{-L^2/2}$ for a standard normal $Z$. Together these give
\[
\operatorname{Re}I_{\mathcal B}\ge\frac{1}{2\pi\sqrt{VW}}\Bigl(e^{-\beta^2/(2W)}-2e^{-L^2/2}-1.65\,\varepsilon\Bigr).
\]

\emph{Minor arcs.} By Lemmas~\ref{lem:K}(i) and~\ref{lem:minor}, $|I_{\mathrm{minor}}|\le4A^2\bigl(e^{-0.4367L^2}+e^{-E_{\mathrm{far}}}\bigr)$.

\emph{Conclusion.} By \eqref{eq:trivial} and Lemma~\ref{lem:saddle}, $2\pi\sqrt{VW}\le16.33m^2$ and $4A^2\le4(s+1)^2\le20.26m^4$. Therefore $I>0$ provided that
\begin{equation}\label{eq:final}
e^{-\beta^2/(2W)}-2e^{-L^2/2}-1.65\,\varepsilon-331\,m^6\bigl(e^{-0.4367L^2}+e^{-E_{\mathrm{far}}}\bigr)>0 .
\end{equation}
We check \eqref{eq:final} for $m\ge10^{30}$.
\begin{itemize}
\item By Lemma~\ref{lem:K}(ii) and Lemma~\ref{lem:shape}, $\beta^2/(2W)\le4s^2/(1.364s^3)\le3.11m^{-1/2}$ in regime (I), and $\beta^2/(2W)\le25m^2/(1.3\cdot10^{-4}m^3)\le1.93\cdot10^5m^{-1}$ in regime (II).
\item By the choice of $L$, $331m^6e^{-0.4367L^2}\le331e^{-10.48}m^{-0.026}\le0.0093$, and $2e^{-L^2/2}\le e^{-500}$.
\item $E_{\mathrm{far}}\ge2.15\cdot10^{-6}\sqrt m\ge2\cdot10^{9}$, so $331m^6e^{-E_{\mathrm{far}}}$ is negligible.
\item Using $s\ge0.9434\sqrt m$, in regime (I) we get $\varepsilon\le351.1L^3m^{-1/4}+3.12L^2m^{-1/2}+R_{\mathrm{tail}}$. At $m=10^{30}$ this is at most $0.34$, since then $L^2\le977.3$ and $R_{\mathrm{tail}}\le e^{-10^5}$ (note that $\sqrt s/L\ge9\cdot10^5$). In regime (II), $\varepsilon\le0.005$.
\end{itemize}
Each of these quantities is decreasing in $m$ for $m\ge10^{30}$, because $\frac{d}{dx}\bigl[\frac32\log(24+13.8x)-\frac x4\bigr]<0$ for $x=\log m>5$. Hence the left-hand side of \eqref{eq:final} is at least $1-10^{-14}-10^{-211}-1.65(0.34)-0.0093-10^{-6}>0.42$ for all $m\ge10^{30}$. This proves \eqref{eq:goal} and the theorem.\qed

\begin{remark}\label{rem:numerics}
The argument shows, more precisely, that uniformly for $5m<t\le T$,
\[
p_t=\frac{F(\rho,r)\rho^{-m}r^{-t'}K(r)\,e^{-\beta^2/(2W)}}{2\pi\sqrt{VW}}\bigl(1+O(L^3m^{-1/4})\bigr).
\]
The threshold $10^{30}$ is very far from the truth. We compared this main term with the exact coefficients computed in Theorem~\ref{thm:computation}. The largest value of $|\log(p_t/\text{main term})|$ over $m<t\le T$ (computed without the factor $e^{-\beta^2/(2W)}$, which is $1+O(1/s)$) is $0.64$, $0.36$, $0.24$, $0.18$, $0.13$ for $m=20,40,80,160,300$. In each case it is attained at $t=m+1$, and it decays roughly like $m^{-1/2}$. Near $t=T$ it is about $1.3/m$. Moreover, the constants in Lemma~\ref{lem:shape} are far from sharp. For example, in regime (I) the true values are $V\approx s$, $W\approx\frac{\pi^2}{3}s^3$ and $M_3\approx9\zeta(3)s^4$, against our bounds $0.3628s$, $0.682s^3$ and $193s^4$. The threshold is driven by the bound $\varepsilon\lesssim351L^3m^{-1/4}$ in regime (I). Sharper moment estimates, an Edgeworth correction exploiting the vanishing of the odd cubic term, and an extension of Corollary~\ref{cor:tail} beyond $t=5m$ should each lower it by many orders of magnitude.
\end{remark}

\section{Graph binomial coefficients and $q$-binomial coefficients}\label{sec:graphs}

Let $G=(V,E)$ be a labelled graph with $n$ vertices. In \cite{EMM}, the number of $k$-combinations of the unlabelled graph $[G]$, i.e.\ the number of orbits of $k$-subsets of $V$ under $\mathrm{Aut}(G)$, was denoted by $\gbin{[G]}{k}$ and called the \emph{graph binomial coefficient}. By Burnside's lemma \cite[Theorem 2.1]{EMM},
\[
\gbin{[G]}{k}=\frac{1}{|\mathrm{Aut}(G)|}\sum_{\sigma\in\mathrm{Aut}(G)}|V^k_\sigma|,
\]
where $V^k_\sigma$ is the set of $k$-subsets fixed by $\sigma$. Explicit formulas were obtained for the path $P_n$ \cite[Corollary 2.1]{EMM} and for the directed cycle $\vec C_n$ \cite[Corollary 2.2]{EMM}:
\begin{equation}\label{eq:EMM}
\gbin{[P_n]}{k}=\begin{cases}\frac12\binom nk & n \text{ even},\ k\text{ odd},\\[2pt] \frac12\binom nk+\frac12\binom{\lfloor n/2\rfloor}{\lfloor k/2\rfloor}&\text{otherwise,}\end{cases}
\qquad
\gbin{[\vec C_n]}{k}=\frac1n\sum_{d\mid \gcd(n,k)}\varphi(d)\binom{n/d}{k/d}.
\end{equation}

The next proposition shows that these two graph binomial coefficients are hidden inside the Gaussian polynomial $\qbin nk$.

\begin{proposition}\label{prop:filters}
Let $0\le k\le n$. Then
\[
\gbin{[P_n]}{k}=\sum_{t\equiv 0\ (\mathrm{mod}\ 2)}[q^t]\qbin nk,
\qquad
\gbin{[\vec C_n]}{k}=\sum_{t\equiv 0\ (\mathrm{mod}\ n)}[q^t]\qbin nk .
\]
\end{proposition}

\begin{proof}
For a polynomial $f$ and $r\ge 1$, the roots-of-unity filter gives
\[\sum_{t\equiv 0\ (\mathrm{mod}\ r)}[q^t]f=\frac1r\sum_{j=0}^{r-1}f(\omega^j),\qquad \omega=e^{2\pi i/r}.\] Let $\zeta$ be a primitive $d$-th root of unity with $d\mid n$. Then
\[
\qbinz nk_{q=\zeta}=\begin{cases}\binom{n/d}{k/d}&\text{if } d\mid k,\\ 0&\text{otherwise}\end{cases}
\]
(see, e.g., \cite[Proposition 4.2]{RSW}). For $r=n$, grouping the $j$'s according to the order $d$ of $\omega^j$ gives $\frac1n\sum_{d\mid \gcd(n,k)}\varphi(d)\binom{n/d}{k/d}$, which is the formula \eqref{eq:EMM} for $\vec C_n$. For $r=2$ we get $\frac12\big(\binom nk+\qbinz nk_{q=-1}\big)$. Moreover $\qbinz nk_{q=-1}=\binom{\lfloor n/2\rfloor}{\lfloor k/2\rfloor}$, unless $n$ is even and $k$ is odd, in which case it is $0$ (this is the case $d=2$ of the evaluation above when $n$ is even; for $n$ odd it follows from $\qbin nk=\qbin{n-1}{k}+q^{n-k}\qbin{n-1}{k-1}$). This gives \eqref{eq:EMM} for $P_n$.
\end{proof}

Proposition~\ref{prop:filters} says that the Burnside averages of \cite{EMM} for $P_n$ and $\vec C_n$ are instances of the $q=-1$ phenomenon of Stembridge \cite{St} and of the cyclic sieving phenomenon of Reiner, Stanton and White \cite{RSW}. Combinatorially, the weight $w(S)=\sum_{s\in S}s-\binom k2$ on $k$-subsets of $\{0,\dots,n-1\}$ is the bridge: a rotation of $\vec C_n$ by $r$ steps changes $w$ by $rk$ modulo $n$, and the reflection of $P_n$ sends $w$ to $k(n-k)-w$.

We now use this point of view for the polynomial $P^{(m)}$. Put $a=2m-1$, $b=3m-1$, $c=m+2$, so that $a+b+c=6m$, and consider the labelled graph
\[
\Gamma_m=P_{c}\ \sqcup\ \vec C_{a}\ \sqcup\ \vec C_{b}
\]
whose vertices are labelled $0,1,\dots,6m-1$ in this order: the path occupies $\{0,\dots,c-1\}$, the first cycle $\{c,\dots,c+a-1\}$ and the second cycle the remaining $b$ labels. The group $\mathbb Z_a\times\mathbb Z_b$ of independent rotations of the two cycles is a subgroup of $\mathrm{Aut}(\Gamma_m)$. For an $m$-subset $S$ of the vertices let
\[
(l,i,j)=(|S\cap P_c|,\ |S\cap \vec C_a|,\ |S\cap\vec C_b|).
\]
By the $q$-Vandermonde convolution, applied twice,
\begin{equation}\label{eq:vandermonde}
\qbin{6m}{m}=\sum_{l+i+j=m}q^{e(l,i,j)}\qbin{c}{l}\qbin{a}{i}\qbin{b}{j},\qquad e(l,i,j)=(c-l)(i+j)+(a-i)j,
\end{equation}
where the $(l,i,j)$-term is the generating function, with respect to $w$, of the $m$-subsets $S$ of type $(l,i,j)$.

Call the type $(l,i,j)$ \emph{free} if $\gcd(i,a)=\gcd(j,b)=1$. If $S$ has free type, then no non-trivial rotation of $\vec C_a$ fixes $S\cap\vec C_a$: a subset of $\mathbb Z_a$ fixed by a rotation of order $d>1$ is a union of cycles of length $d$, so $d\mid i$ and $d\mid a$. Similarly for $\vec C_b$. Hence the $\mathbb Z_a\times\mathbb Z_b$-orbit of $S$ has exactly $ab$ elements.

\begin{theorem}\label{thm:orbit}
Let $m\ge 2$ with $m\not\equiv 2,3\pmod 5$, and define
\[
F^{(m)}(q)=\sum_{\substack{l+i+j=m\\ \gcd(i,a)=\gcd(j,b)=1}}q^{e(l,i,j)}\qbin{c}{l}\cdot\frac{1}{[a]_q}\qbin{a}{i}\cdot\frac{1}{[b]_q}\qbin{b}{j},
\qquad E^{(m)}=P^{(m)}-F^{(m)}.
\]
Then:
\begin{enumerate}
\item $F^{(m)}(q)$ is a polynomial with non-negative integer coefficients, and $F^{(m)}(1)$ is the number of $\mathbb Z_a\times\mathbb Z_b$-orbits of $m$-combinations of $\Gamma_m$ of free type;
\item $E^{(m)}(q)$ is a polynomial, and
\[
[a]_q[b]_q\,E^{(m)}(q)=\sum_{\substack{l+i+j=m\\ \gcd(i,a)>1\text{ or }\gcd(j,b)>1}}q^{e(l,i,j)}\qbin{c}{l}\qbin{a}{i}\qbin{b}{j};
\]
\item consequently, Conjecture~\ref{conj:general} for $m$ is implied by the statement that $E^{(m)}(q)+q+q^{D_m-1}$ has non-negative coefficients.
\end{enumerate}
\end{theorem}

\begin{proof}
(1) If $\gcd(i,a)=1$ then $\gcd(i,a-i)=1$, and by \cite[Lemma 5.2]{GK} the quotient $\frac{1-q}{1-q^{a}}\qbin{a}{i}=\frac1{[a]_q}\qbin ai$ is a polynomial with non-negative coefficients (it is the rational $q$-Catalan number of type $(i,a-i)$). The same holds for $b,j$. Hence each summand of $F^{(m)}$ is in $\N[q]$. At $q=1$ each summand equals $\binom cl\binom ai\binom bj/(ab)$, which is the number of free-type subsets of type $(l,i,j)$ divided by the common orbit size $ab$.

(2) This follows from \eqref{eq:vandermonde}, Proposition~\ref{prop:poly} and part (1), since $E^{(m)}$ is the difference of two polynomials.

(3) This is immediate from $P^{(m)}=F^{(m)}+E^{(m)}$ and $F^{(m)}\ge 0$.
\end{proof}

Thus the free orbits, which form the ``generic'' part of the Burnside count for $\Gamma_m$, contribute only non-negative coefficients. All possible negativity is concentrated in the $m$-combinations whose intersection with one of the two cycles has a size sharing a factor with the length of that cycle. In particular, this includes the combinations that avoid one of the cycles completely, which are fixed by all rotations of that cycle. When $a$ and $b$ are both prime (for example $m=10,16$), these are exactly the combinations with $i=0$ or $j=0$.

\begin{proposition}\label{prop:Edata}
For every $m\in\{9,10,11,14,15,16,19,20,21,24,25,26,29,30\}$, the polynomial $E^{(m)}$ has exactly two negative coefficients, namely $-1$ at $q^1$ and at $q^{D_m-1}$.
\end{proposition}

\begin{proof}
Direct computation, as in Theorem~\ref{thm:computation}.
\end{proof}

\begin{table}[h]
\centering
\begin{tabular}{rrr}
\toprule
$m$ & $F^{(m)}(1)$ & $E^{(m)}(1)$\\
\midrule
$10$ & $134\,716\,109$ & $2\,115\,157$\\
$16$ & $454\,084\,207\,510\,368$ & $447\,078\,960\,834$\\
$20$ & $8\,522\,803\,227\,830\,696\,313$ & $4\,281\,293\,813\,097\,872\,013$\\
$24$ & $\approx 4.007\cdot 10^{23}$ & $\approx 1.040\cdot 10^{19}$\\
$30$ & $\approx 2.524\cdot 10^{30}$ & $\approx 4.299\cdot 10^{24}$\\
\bottomrule
\end{tabular}
\caption{Sizes of the free part and of the correction term.}
\label{tab:FE}
\end{table}

When $a$ and $b$ are both prime ($m=10,16,24,30$ in Table~\ref{tab:FE}), the free part dominates. When $a$ or $b$ has small prime factors, more types fail to be free and $E^{(m)}$ becomes comparable to $F^{(m)}$ (for example $m=25$, where $a=49$ and $b=74$).

\begin{remark}
A finer decomposition uses all free orbits, not only those of free type. By Möbius inversion, $\qbin ai=\sum_{d\mid\gcd(a,i)}\mathrm{Fr}_{a/d,i/d}(q^d)$, where
\[
\mathrm{Fr}_{a,i}(q)=\sum_{d\mid \gcd(a,i)}\mu(d)\qbinz{a/d}{i/d}_{q^d}
\]
is a $q$-analogue of the number of aperiodic $i$-subsets of $\mathbb Z_a$. We checked that $\mathrm{Fr}_{n,i}(q)/[n]_q\in\N[q]$ for all $1\le i<n\le 30$. Replacing $\qbin ai$ and $\qbin bj$ by $\mathrm{Fr}_{a,i}$ and $\mathrm{Fr}_{b,j}$ in $F^{(m)}$ gives a non-negative part whose complement is very small. For instance, the ratio of the complement to $P^{(m)}(1)$ is about $1.6\%$, $0.17\%$ and $0.02\%$ for $m=10,15,20$, and this complement again has exactly the two coefficients $-1$ for $m\in\{9,10,11,14,15,16,19,20\}$. It would be interesting to prove that $\mathrm{Fr}_{n,i}(q)/[n]_q\in\N[q]$ in general and to give it a combinatorial interpretation.
\end{remark}

\section{Concluding remarks}

The phenomenon described in Conjecture~\ref{conj:general} does not seem to be special to the parameters $(6m,m,2m-1,3m-1)$. A computer search over all $N\le 45$ shows that, among the $10\,053$ quotients $\qbin Nk/([a]_q[b]_q)$ with $2\le k\le N/2$ and $2\le a\le b\le N$ that are polynomials, $8\,582$ have exactly the two coefficients $-1$ at $q^1$ and at $q^{D-1}$ and no other negative coefficient. The exceptions occur mainly when $a$ is small compared with $k$; in that case the coefficient $-1$ of $q^1$ is propagated to $q^{1+a},q^{1+2a},\dots$ and is not compensated. For instance, for the families
\[
\frac{1}{[2m-1]_q[6m-1]_q}\qbin{12m}{m},\qquad
\frac{1}{[4m-1]_q[6m-1]_q}\qbin{12m}{5m},\qquad
\frac{1}{[4m-1]_q[6m-1]_q}\qbin{12m}{3m},
\]
every case with $2\le m\le 10$ in which the quotient is a polynomial shows exactly the Guo--Krattenthaler pattern. Lemma~\ref{lem:defect} applies verbatim to such quotients.

The most pressing open point is the finite gap left by Theorem~\ref{thm:asymptotic}.

\begin{problem}
Lower the threshold $10^{30}$ of Theorem~\ref{thm:asymptotic} far enough to meet the range $m\le300$ covered by Theorem~\ref{thm:computation}. The threshold is driven by regime (I) of Section~\ref{sec:asymptotic}, i.e.\ by $t$ just above $5m$, where the effective number of active Bernoulli variables is only of order $\sqrt t$, and by the crude constants of Lemma~\ref{lem:shape}. Corollary~\ref{cor:tail} disposes of the range $t\le 5m$ by a purely combinatorial argument. Extending it to $t\le\varepsilon m^2$ would mean controlling the ``box corrections'' $\prod_{j=1}^m(1-q^{5m+j})$ in $\qbin{6m}{m}=\prod_{j=1}^{m}(1-q^{5m+j})/(q;q)_m$. Such an extension would leave only the moderate- and central-deviation ranges to the analytic argument, where the effective number of variables is of order $m$.
\end{problem}

\begin{problem}
Characterize the quadruples $(N,k,a,b)$ for which $\qbin Nk/([a]_q[b]_q)$ is a polynomial whose only negative coefficients are $-1$ at $q^1$ and at $q^{D-1}$. In particular, is this always the case when $\qbin Nk/[a]_q$ is weakly increasing from $q^1$ to its middle coefficient?
\end{problem}

By Lemma~\ref{lem:defect}, the answer to the last question is affirmative as soon as, in addition, the coefficients of $q^{b}$ and $q^{b+1}$ in $\qbin Nk/[a]_q$ are distinct.

\end{document}